\documentclass[a4paper,10pt,reqno]{amsart}
\usepackage[centertags]{amsmath}
\usepackage{amsfonts,amssymb,amsthm}          
\usepackage{newlfont}
\usepackage[english]{babel} 
\usepackage[body={15cm, 23cm},centering]{geometry} 
\usepackage{fancyhdr}
 \theoremstyle{plain}
 \begingroup
 \newtheorem{teo}{Theorem}[section]
 \newtheorem{lem}[teo]{Lemma}
 \newtheorem{prop}[teo]{Proposition}
 
 \endgroup 
  
 \theoremstyle{definition}
 \begingroup
 \newtheorem{defi}[teo]{Definition}
 \newtheorem{ex}[teo]{Example}
 \newtheorem{ass}{Assumption}
 \endgroup 
  
 \theoremstyle{remark}
 \begingroup
 \newtheorem{rem}[teo]{Remark}
 \endgroup

 \numberwithin{equation}{section}

\renewcommand{\b}{\mathcal{B}}
\newcommand{\f}{\mathcal{F}}
\newcommand{\m}{\mathcal{M}} 
\newcommand{\p}{\mathcal{P}} 
\newcommand{\leb}{\mathcal{L}} 

\newcommand{\jm}{j}

\newcommand{\e}{\varepsilon}
\newcommand{\Om}{\Omega}

\newcommand{\R}{{\mathbb R}}

\newcommand{\NP}{\mathbb{N}^{+}}

\newcommand{\M}{\mathbb{M}}

\newcommand{\Y}{\mathcal{Y}}
\newcommand{\YD}{\widehat{\mathcal{Y}}}
\newcommand{\y}{y^1, \ldots,y^n}
\newcommand{\yk}{y^1, \ldots,y^k}

\newcommand{\g}{\nabla}

\newcommand{\norm}[1]{\left\Vert#1\right\Vert}
\newcommand{\snorm}[1]{\left|#1\right|}

\DeclareMathOperator{\essinf}{ess\,inf}

\newcommand{\ci}{C_c^\infty\negthinspace\left(\Om\right)}
\newcommand{\cip}{C_{per}^\infty(\square^n)}

\title[Homogenization of Functionals with Discontinuous Integrand]
{Multiscale Homogenization of Convex Functionals with Discontinuous Integrand}
\author[Marco Barchiesi]{Marco Barchiesi}
\address[Marco Barchiesi]{S.I.S.S.A., Via Beirut 2-4, 34014, Trieste, Italy}
\email[Marco Barchiesi]{barchies@sissa.it}

\begin{document}
\maketitle

\begin{center}
\begin{minipage}{12cm}
\small{
\noindent {\bf Abstract.} 
This article is devoted to obtain the $\Gamma$-limit, as $\e$ tends to zero, of the family of functionals 
\begin{equation*}
u\mapsto\int_{\Om}f\Bigl(x,\frac{x}{\e}, \ldots, \frac{x}{\e^n}, \g u(x)\Bigr)dx ,
\end{equation*}
where $f=f(x, \y, z)$ is periodic in $\y$, convex in $z$ and satisfies
a very weak regularity assumption with respect to $x,\y$.
We approach the problem using the multiscale Young measures. 

\vspace{10pt}

\vspace{6pt}
}
\end{minipage}
\end{center}

\bigskip
\tableofcontents

\section{Introduction}

\noindent
Multiscale composites are structures constituted by two or more materials which are finely mixed
on many different microscopic scales.
The fact that a composite often combines the properties of the constituent materials
makes these structures particularly interesting in many fields of science.
There is a vast literature on the subject; we refer the reader 
to \cite{Mil02} and references therein.

Determining macroscopic behavior of these strongly heterogeneous structures when the size $\e$ of the 
heterogeneity becomes ``small'' is the aim of homogenization theory.

In the particular case of a periodic multiscale composite,
from a variational point of view, the homogenization problem  
is to characterize the behavior, for the parameter $\e$ tending to zero,
of functionals on $W^{1,p}(\Om, \R^m)$ of the type
\begin{equation}\label{prima}
F_\e (u)=\int_\Om f\Bigl(x,\langle\frac{x}{\rho_1(\e)}\rangle, \ldots, \langle\frac{x}{\rho_n(\e)}\rangle, \g u(x)\Bigr)dx\,,
\end{equation}

\vspace{8pt}
\noindent
where $\langle\cdot\rangle$ denotes the fractional part of a vector componentwise,
$\Om$ is an open bounded domain in $\R^d$, $\square$ is the unit cell $[0, 1)^d$, $\rho_k$ are the length scales
and $f=f\left(x, \y, z\right)$ is a non-negative function on $\Om\times\square^n\times\M^{m\times d}$. 

The purpose of this paper is to analyze \eqref{prima} under the following assumptions.
\begin{ass}\label{ass1} $f$ is convex in the argument $z$ for all $x\in\Om$ and $\y\in\square$.
\end{ass}
\begin{ass}\label{ass2} $f$ is $p$-coercive and with $p$-growth:
\begin{equation*}
c_1\snorm{z}^p\leq f\left(x,\y, z\right)\leq c_2\left(1+\snorm{z}^p\right)
\end{equation*}
for some $p\in(1,+\infty)$, \,$c_1, c_2>0$ and for all $\left(x, \y, z\right)\in\Om\times\square^n\times\M^{m\times d}$.
\end{ass}
\begin{ass}\label{ass3} $f$ is an \emph{admissible integrand}, \textit{i.e.},
for every $\delta>0$ there exist a compact set $X\subseteq\Om$ with $\snorm{\Om\backslash X}\leq\delta$ and 
a compact set $Y\subseteq\square$ with $\snorm{\square\backslash Y}\leq\delta$, 
such that $f\vert _{X\times Y^n\times\M^{m\times d}}$ is continuous.
\end{ass}

In particular we cover the following two significant cases (see Examples \ref{int ammiss} and \ref{int ammiss due}).
\begin{itemize}
\item[(i)] The case of a single microscale ($n=1$):
the function $f:\Om\times\square\times\M^{m\times d}\rightarrow[0,+\infty)$ 
is continuous in $x$, measurable in $y$ and satisfies Assumptions 
\ref{ass1} and \ref{ass2}. Notice that $f$ is continuous in~$z$ uniformly with respect to~$x$ and 
hence is continuous in $(x,z)$. It is possible to interchange the regularity conditions on $f$
requiring the measurability in $x$ and the continuity in $y$.
\item[(ii)] The case of a multiscale mixture of two materials:
the function $f:\Om\times\square^n\times\M^{m\times d}\rightarrow[0,+\infty)$ is  of the type
\begin{equation}\begin{split}\label{composito}
f(x,\y,z)=\prod_{k=1}^{n} \chi_{P_k}(y^k)\,f_1(x,\y,z)&\\+\Bigl[1-\prod_{k=1}^{n} \chi_{P_k}(y^k)\Bigr]&\,f_2(x,\y,z)\,,
\end{split}\end{equation}
where $\chi_{P_k}$ ($k=1,\ldots,n$) is the characteristic function of a measurable subset $P_k$ of \,$\square$
and the functions $f_1,f_2:\Om\times\square^n\times\M^{m\times d}\rightarrow[0,+\infty)$ are 
measurable in $x$, continuous in $(y^1,\ldots,y^n)$ and satisfy Assumptions \ref{ass1} and \ref{ass2}.
The regularity conditions on $f_1$ and $f_2$ can be replaced by the continuity in $(x,y^1,\ldots,y^{n-1})$
and the measurability in the fastest oscillating variable $y^n$. 
\end{itemize}

\noindent
Problems of the type \eqref{prima} have captured the attention of many authors.
For instance, the case of a single microscale 
\begin{equation*}
\int_\Om f\Bigl(x, \langle\frac{x}{\e}\rangle, \g u(x)\Bigr)dx
\end{equation*}
has been studied by Braides (see \cite{Br83} and also \cite[Chapter 14]{Br98})
under Assumption \ref{ass1} and requiring in addition a $p$-growth condition on the integrand $f$ 
and a uniform continuity in $x$, precisely
\begin{equation}\label{uniform}
\snorm{f(x,y,z)-f(x',y,z)}\leq\omega(\snorm{x-x'})\bigl[\alpha(y)+f(x,y,z)\bigr]
\end{equation}
for all $x,x'\in\R^d$, $y\in\square$ and $z\in\M^{m\times d}$,
where $\alpha\in L^1(\square)$ and $\omega$ is a continuous positive function with $\omega(0)=0$.
Recently Ba{\'\i}a and Fonseca \cite{bai05} have studied this problem under Assumpion
\ref{ass2} and requiring continuity in $(y,z)$ and measurability in $x$.

In \cite{Br00bis} (see also \cite[Chapter 22]{Br98}, \cite{Lukk99} and \cite{Mil00}) Braides and Lukkassen 
study functionals of the form 
\begin{equation*}
\int_\Om f\Bigl(\langle\frac{x}{\e}\rangle, \ldots,\langle\frac{x}{\e^n}\rangle, \g u(x)\Bigr)dx\,.
\end{equation*} 
The authors provide an iterated homogenization formula for functions as in \eqref{composito}
with an additional request on the functions $f_1$ and $f_2$ of a uniform continuity, similar to \eqref{uniform}, 
with respect to the slower oscillating variables $y^1, \ldots,y^{n-1}$. 
The same result is obtained by Fonseca and Zappale \cite{Fons03}
but with a continuous function $f$ satisfying Assumptions \ref{ass1} and \ref{ass2}.

Since the variable $x$ describes the macroscopic heterogeneity of the constituent materials while 
the variables $\y$ describe the microscopic heterogeneity of the composite structure, 
it is desirable  to have the weakest possible regularity on them.
In particular, the oscillating variables should be able to describe the discontinuity on the interfaces
between different materials.
At any rate the only request that $f$ is borelian is not enough to obtain a homogenization formula, 
as it is shown  in Examples \ref{final uno} and \ref{final due} (see also \cite{All92} and \cite{Don83}).

In order to weaken the continuity assumptions taken in the works cited above,
 we approach the problem using the multiscale Young measures 
as in \cite{Ped05} (see also \cite{Mas01} and \cite{Ped04}).
The peculiarity of our work is the introduction of the concept of \emph{admissible integrand} (Definition \ref{new}).
The crucial point is to extend the lower semicontinuity property \eqref{lower sem} to this kind of integrand:
this is achieved in Theorem \ref{continuita2}.

The paper is organized as follows. In Section 2 we recall concepts and basic facts about Young measures.
In Section 3 we introduce the notion of multiscale convergence in the general framework of multiscale Young measures.
In Section 4 we discuss the properties of admissible integrands.
By Theorems \ref{continuita1} and \ref{continuita2} we derive, in Section~5, the upper and
lower estimates for the $\Gamma(L^p)$-limit of the family $F_\e$ (Lemmas \ref{sup lemma} and \ref{inf lemma}). 
Finally, in Section 6, we give an iterated homogenization formula.

\section{Young measures}
\noindent
We gather briefly in this section some of the main results about Young measures, for more details and 
proofs we refer the reader to \cite{B88} and \cite{Val90}. 

We denote with
\begin{itemize}
\item $D$ a bounded Lebesgue measurable subset of $\R^l$ ($l\geq1$), equipped with the Lebesgue \\$\sigma$-algebra $\f(D)$;
\item $\snorm{A}$ the Lebesgue measure of a set $A\in\f(D)$;
\item $S$ a locally compact, complete and separable metric space, equipped with \\the Borel $\sigma$-algebra $\b(S)$;
\item $\leb(D, S)$ the family of measurable functions $u:D\rightarrow S$;
\item $C_0(S)$ the space $\lbrace\phi: S\rightarrow\R \;\text{continuous} : \forall\delta>0 \;\exists K\subseteq S
\;\text{compact} :\\ \vert\phi(z)\vert<\delta \;\text{for} \;z\in S\setminus K\rbrace$, endowed with the supremum norm;
\item $\m(S)$ the space of Radon measures on $S$;
\item $\p(S):=\lbrace\mu\in\m(S) : \mu\geq0 \;\text{and} \;\mu(S)=1\rbrace$ the set of probability measures on $S$;
\item $L^1\negthinspace\left(D, C_0(S)\right)$ the Banach space of all measurable maps 
 $x\in D\xrightarrow{\phi} \phi_x\in C_0(S)$ such that the quantity 
 $\norm{\phi}_{L^1}:=\int_D\norm{\phi_x}_{C_0(S)}dx$
 is finite;
\item $L^{\infty}_w\left(D, \m(S)\right)$ the Banach space of all weak* measurable maps 
 $x\in D\xrightarrow{\mu}\mu_x\in\m(S)$ such that  
 $\norm{\mu}_{L^{\infty}_w}:=\operatorname{ess}\sup_{x\in D}\norm{\mu_x}_{\m(S)}$
 is finite;
\item $\Y(D, S)$ the family of all weak* measurable maps 
 $\mu:D\negthinspace\rightarrow\negthinspace\m(S)$ such that $\mu_x\in\p(S)$ \textit{a.e.} $x\in D$.
 \end{itemize}

\begin{rem}
\begin{itemize}\label{oss princ}
\item[]
\item[(i)]As it is known, the dual of $C_0(S)$ may be identified with the set of $S$-valued Radon
 measures through the duality
\begin{equation*}
\left\langle \mu,\phi\right\rangle=\int_{S}\phi\,d\mu \qquad\forall \mu\in\m(S) 
\quad\mathrm{and} \quad\forall \phi\in C_0(S)\,.
\end{equation*}
\item[(ii)]A map $\mu: D\rightarrow\m(S)$ is said to be \emph{weak* measurable} if 
    $x\rightarrow\left\langle \mu_x,\phi\right\rangle$ is measurable for all $\phi\in C_0(S)$.
     In particular $x\rightarrow\norm{\mu_x}_{\m(S)}$ is measurable.
\item[(iii)]More precisely, the elements of $L^1\negthinspace\left(D, C_0(S)\right)$,
     $L^{\infty}_w\left(D, \m(S)\right)$ and $\Y(D, S)$ are equivalence
     classes of maps that agree \textit{a.e.}; we usually do not distinguish 
     these maps from their equivalence classes.
\item[(iv)]$L^{\infty}_w\left(D, \m(S)\right)$ can be identified with the dual 
    of $L^1\negthinspace\left(D, C_0(S)\right)$ through the duality
\begin{equation*}
\left\langle \mu,\phi\right\rangle=\int_D\left\langle \mu_x,\phi_x\right\rangle dx 
 \qquad \forall \mu\in L^{\infty}_w\bigl(D, \m(S)\bigr) \quad\mathrm{and} 
 \quad\forall\,\phi\in L^1\negthinspace\bigl(D, C_0(S)\bigr).
\end{equation*}
In the following we will refer to the weak* topology of $L^{\infty}_w\left(D, \m(S)\right)$ as the topology 
induced by this duality pairing.
\item[(v)] Let $\YD(D, S):=\lbrace\varrho\in\m(D\times S): \varrho\geq0 
\;\text{and} \;\varrho(A\times S)=\snorm{A} \;\forall A\in\b(D)\rbrace$. 
By the Disintegration Theorem \cite{Val73}, the map
which associates to $\mu\in\Y(D, S)$ the measure $\widehat{\mu}\in\YD(D, S)$
defined by
\begin{equation*}
\widehat{\mu}(A):=\int_D\left(\int_S\chi_A(x, z)d\mu_x(z)\right)dx  \quad\forall A\in\b(D\times S)
\end{equation*}
induces a bijection between $\Y(D, S)$ and $\YD(D, S)$.
Given a function $f\negthinspace:D\times S\rightarrow\R$ \,$\widehat{\mu}$-integrable, it turns out that
$f(x,\cdot)$ is $\mu_x$-integrable for \textit{a.e.} $x\in\negthinspace D$, 
$x\rightarrow\int_Sf(x, z)d\mu_x(z)$ is integrable and
\begin{equation*}
\int_{D\times S}fd\widehat{\mu}=\int_D\left(\int_S f(x, z)d\mu_x(z)\right)dx\,;
\end{equation*}

\vspace{4pt}
\noindent
this last equality remains true if $f$ is $\f(D)\otimes\b(S)$-measurable and non-negative.
\end{itemize}
\end{rem}

\vspace{10pt}
\noindent
The family $\leb(D, S)$ can be embedded in $L^{\infty}_w\left(D, \m(S)\right)$ associating
to every $u\in\leb(D, S)$ the function 
\begin{equation*}
x\xrightarrow{\delta_u}\delta_{u(x)} ,
\end{equation*} 
where $\delta_{u(x)}$ is the Dirac probability measure concentrated at the point $u(x)$.
\begin{defi}
A function $\mu\in L^{\infty}_w\left(D, \m(S)\right)$ is called the \emph{Young measure} generated by
the sequence $u_h$ if $\delta_{u_h}\rightharpoonup\mu$ in the weak* topology.
\end{defi}
\begin{rem}
This notion makes sense: by the identification of 
$L^{\infty}_w\left(D, \m(S)\right)\simeq L^1\left(D, C_0(S)\right)^*$
and as a direct consequence of the Banach-Alaoglu theorem, every sequence $u_h$ in $\leb(D, S)$ admits
a subsequence generating a Young measure.
\end{rem}

The following result is a ``light" version of the Fundamental Theorem on Young Measures.
\begin{teo}\label{fund young}
Let $u_h$ be a sequence in $\leb(D, S)$ generating a Young measure $\mu$ and 
for which the ``tightness condition'' is satisfied, \textit{i.e.}, 
\begin{equation}\label{tight}
\forall\delta>0 
\;\exists K_\delta\subseteq \,S\,\text{compact}:  
\sup_{h\in\NP}\bigl\vert\lbrace x : u_h(x)\not\in K_\delta\rbrace\bigr\vert\leq\delta\,.
\end{equation}
The following properties hold:
\begin{enumerate}
\item[i)]$\mu\in\Y(D, S)$;
\item[ii)] if $f:D\times S\rightarrow[0, +\infty)$ is a Caratheodory integrand, then
\begin{equation*}
\liminf_{h\rightarrow+\infty}\int_\Om f\bigl(x, u_{h}(x)\bigr)dx\geq\int_\Om\overline{f}(x)dx
\end{equation*}
where
\begin{equation*}
\overline{f}(x):=\int_Sf(x, z)d\mu_x(z);
\end{equation*}
\item[iii)] if $f:D\times S\rightarrow\R$ is a Caratheodory integrand and
$f\left(\cdot, u_{h}(\cdot)\right)$ is equi-integrable, then
 $f$ is $\widehat{\mu}$-integrable and
$f\left(\cdot, u_{h}(\cdot)\right)\rightharpoonup\overline{f}$ weakly in $L^1(\Om)$.
\end{enumerate}
\end{teo}

\noindent
We remember that a $\f(D)\otimes\b(S)$-measurable function $f$ is a Caratheodory integrand if
$f(x,\cdot)$ is continuous for all $x\in D$.

\section{Multiscale Young measures}

\noindent
We introduce now the notion of multiscale convergence, an extension of the two-scale convergence
carried out by Allaire (\cite{All96}) in joint work with Briane. We present it in the general
framework of multiscale Young measures, following essentially the ideas exposed 
in \cite{Val97}, \cite{Bald99} and \cite{Mas01}.

We start presenting an example that does not only show a fine and explicit case of Young measure,
but it is a fundamental mainstay in this section.
Before we add some new notations:
\begin{itemize}
\item \,$\Omega$ is a bounded open subset of ${\R}^d$ ($d\geq1$), equipped with the Lebesgue $\sigma$-algebra $\f(\Om)$;
\item $\square$ is the unit cell $[0, 1)^d$, equipped with the Lebesgue $\sigma$-algebra $\f(\square)$; 
\item $n$ is the number of scales, a positive integer;
\item $\rho_1, \ldots, \rho_n$ are positive functions of a parameter $\e>0$ which converge to $0$ as $\e$ does,
for which the following \emph{separation of scales hypothesis} is supposed to hold:
\begin{equation*}
\lim_{\e\rightarrow0^+}\frac{\rho_{k+1}(\e)}{\rho_k(\e)}=0 
\quad\forall k\in\{1,\ldots ,n-1\};
\end{equation*}
\item $p\in(1,+\infty)$ and $q\in[1,+\infty]$ (unless otherwise stated), moreover $q'$ is the H\"olderian 
conjugate exponent of $q$;
\item $C^\jm_c(\Om)$ stands for the space of $\jm$-differentiable functions in $\Om$ with compact support;
\item $C^\jm_{per}(\square^k)$ is the space of the functions $u=u(y^1,\ldots, y^k)$ in
$C^\jm((\R^d)^k)$ separately $\square$-periodic in $y^1,\ldots, y^k$;
\item $W^{1,q}_{per}(\square^k)$ 
denotes the space of the functions $u=u(y^1,\ldots, y^k)$ in $W^{1,q}_{loc}((\R^d)^k)$ 
separately $\square$-periodic in $y^1,\ldots, y^k$.
\end{itemize}

\vspace{5pt}
\noindent
We fix a sequence $\e_h\rightarrow0^+$ of values of the parameter $\e$.

\vspace{7pt}
\begin{ex}\label{esempio young}
We denote by $T$ the set \,$\square$ \,equipped with the topological and differential
structure of the $d$-dimensional torus and with the Borel $\sigma$-algebra $\b(T)$;
any function on $T$ can be identified with
its periodic extension to $\R^d$, in particular
\begin{equation*}
C(T)=C_0(T)\simeq C_{per}(\square).
\end{equation*}

\vspace{10pt}
\noindent
We consider the sequence $v_h:\Om\rightarrow\square^n$ defined by
\begin{equation}\label{v}
v_h(x):=\left(\langle\frac{x}{\rho_1(\e_h)}\rangle, \ldots, \langle\frac{x}{\rho_n(\e_h)}\rangle\right)\,.
\end{equation}
Here $\langle\cdot\rangle$ denotes the fractional part of a vector componentwise.
For our example, we need an auxiliary ingredient concerning weak convergence. 
It is a particular case of \cite[Proposition 3.3]{Don83bis}.
\begin{teo}\textnormal{Riemann-Lebesgue lemma:}
given $\phi\in C_{per}(\square^n)$, define $\phi_h(x):=\phi(v_h(x))$. 
Then $\phi_h\rightharpoonup\int_{\square^n}\phi\left(\y\right)dy^1\ldots dy^n$
weakly* in $L^{\infty}(\Om)$.
\end{teo}

\vspace{8pt}
As consequence of Riemann-Lebesgue lemma, for all $\varphi\in L^1(\Om)$ and $\phi\in C_{per}(\square^n)$ 
\begin{equation*}
\int_\Om\varphi(x)\phi\left(v_h(x)\right)dx
\rightarrow
\int_{\Om\times\square^n}\varphi(x)\phi\left(\y\right)dx\,dy^1\ldots dy^n .
\end{equation*}
The map $\varphi\otimes\phi$ that takes every $x\in\Om$ into $\varphi(x)\phi(\cdot)\in C_{per}(\square^n)$ 
belongs to $L^1\left(\Om, C_{per}(\square^n)\right)$.
Since the space
$L^1(\Om)\otimes C_{per}(\square^n)$, defined as the linear closure of 
$\lbrace\varphi\otimes\phi : \varphi\in L^1(\Om) \;\text{and} \;\phi\in C_{per}(\square^n)\rbrace$,
is dense in $L^1\left(\Om, C_{per}(\square^n)\right)$, we conclude that $v_h$ generates the 
Young measure $\mu\in\Y(\Om, T^n)$ with
\begin{equation*}
\mu_x=\leb_{\llcorner\square^n} \;\;\textrm{for \textit{a.e.}} \,\;x\in\Om\,,
\end{equation*}
where $\leb_{\llcorner\square^n}$ is the restriction to $\square^n$ of the Lebesgue
measure on $\left(\R^d\right)^n$.
\end{ex}

\vspace{10pt}
\begin{defi}
Let $u_h$ be a sequence in $L^1(\Om)$. The sequence $u_h$ is said to be multiscale
convergent to a function $u=u(x, y^1, \ldots, y^n)\in L^1(\Om\times\square^n)$ if
\begin{equation*}\begin{split}
\lim_{h\rightarrow+\infty}\int_{\Om}\varphi(x)&\phi\left(\frac{x}{\rho_1(\e_h)}, \ldots,\frac{x}{\rho_n(\e_h)}\right)u_h(x)dx\\
&=\int_{\Om\times\square^n}\varphi(x)\phi\left(\y\right)u\left(x, \y\right)dx\,dy^1\ldots dy^n 
\end{split}\end{equation*}
for any $\varphi\in\ci$ and any $\phi\in\cip$. We simply write $u_h\rightsquigarrow u$.
A sequence in $L^1(\Om, \R^m)$ is called multiscale convergent if it is so componentwise.
\end{defi}
\begin{prop}\label{conv debole}
Let $u_h$ be an equi-integrable sequence in $L^1(\Om)$ multiscale convergent 
to a function $u\in L^1(\Om\times\square^n)$. Then $u_h$ converges weakly to $u_\infty$ in $L^1(\Om)$,
where 
\begin{equation*}
u_\infty(x):=\int_{\square^n}u\left(x, \y\right)dy^1\ldots dy^n\,.
\end{equation*}
\end{prop}
\begin{proof}
An equi-integrable sequence is  sequentially  weakly compact in $L^1$, therefore
it is sufficient to prove that $u_h\rightarrow u_\infty$ in distribution.
But this is a direct consequence of the definition, taking $\phi\equiv1$\,.
\end{proof}

\vspace{6pt}
Let $u_h$ be a bounded sequence in $L^1(\Om, \R^m)$; we consider  
the sequence $w_h:\Om\rightarrow\square^n\times\R^m$ defined by
\begin{equation}\label{w}
w_h(x):=
\left(\langle\frac{x}{\rho_1(\e_h)}\rangle, \ldots, \langle\frac{x}{\rho_n(\e_h)}\rangle, u_h(x)\right)\,.
\end{equation}

\vspace{6pt}
\noindent
Suppose that $w_h$ generates a Young measure $\mu$ (at any rate this is true, up to a subsequence).
Thanks to the boundness hypothesis,
it can be easily proved that $w_h$ satisfies tightness condition~\eqref{tight}, so $\mu\in\Y(\Om, T^n\times\R^m)$.
Roughly speaking, by Remark \ref{oss princ}(v), it is possible to piece together $\mu$ in a measure
$\widehat{\mu}\in\YD(\Om, T^n\times\R^m)$. Thanks to Example \ref{esempio young}, actually
$\widehat{\mu}\in\YD(\Om\times\square^n, \R^m)$ and so, by Remark \ref{oss princ}(v) again,
it is possible to dismantle this measure in a new function $\nu\in\Y(\Om\times\square^n, \R^m)$, called the
\emph{multiscale Young measure} generated by $u_h$. 
In particular we have:
\begin{teo}\label{multi yuong}
Let $\mu\in\Y(\Om, T^n\times\R^m)$ be the Young measure generated by $w_h$ and 
let $\nu\in\Y(\Om\times\square^n, \R^m)$ be the multiscale Young measure generated by $u_h$. Then
\begin{equation*}\begin{split}
\int_\Om\left(\int_{\square^n\times\R^m}f(x, \y, z)d\mu_x(\y, z)\right)&dx\\
&\hspace*{-100pt}=\int_{\Om\times\square^n}\left(\int_{\R^m}f(x, \y, z)d\nu_{(x, \y)}(z)\right)dx\,dy^1\ldots dy^n 
\end{split}\end{equation*}
for all $f:\Om\times\square^n\times\R^m\rightarrow\R$\; \;$\widehat{\mu}$-integrable or non-negative
$\f(\Om)\otimes\b(T^n\negthinspace\times\R^m)$-measurable.
\end{teo}

\vspace{6pt}
The next statement lights up the link between Young measures and multiscale convergence.
Sometimes we will use in the sequel the shorter notation $y:=(\y)$\,.
\begin{teo}\label{fund multi young}
Let $u_h$ be a bounded sequence in $L^q(\Om, \R^m)$, $q\in[1,+\infty)$, generating a multiscale Young measure $\nu$.
The following properties hold:
\begin{enumerate}
\item[i)] the center of mass $\overline{\nu}$, defined by
\begin{equation*}
\overline{\nu}(x, \y):=\int_{\R^m}z\;d\nu_{(x, \y)}(z)\,,
\end{equation*}
is in $L^q(\Om\times\square^n, \R^m)$; 
\item[ii)]if $u_h$ is equi-integrable, then $u_h\rightsquigarrow\overline{\nu}$;
\item[iii)] if $f:\Om\times T^n\times\R^m\rightarrow[0, +\infty)$ is a Caratheodory integrand, \textit{i.e.},
$\f(\Om)\otimes\b(T^n\negthinspace\times\R^m)$-measurable and continuous on $T^n\times\R^m$, then
\begin{equation}\label{lower sem}
\liminf_{h\rightarrow+\infty}\int_{\Om}
f\bigl(x, w_h(x)\bigr)dx\geq\int_{\Om\times\square^n}\overline{f}(x, \y)\,dx\,dy^1\ldots dy^n ,
\end{equation}
where 
\begin{equation*}
\overline{f}(x, \y):=\int_{\R^m}f(x, \y, z)\;d\nu_{(x, \y)}(z)\,;
\end{equation*}
\item[iv)] if $f:\Om\times T^n\times\R^m\rightarrow\R$ is a Caratheodory integrand and
$f(\cdot, w_h(\cdot))$ is equi-integrable, then $f(x, y, \cdot)$ is $\nu_{(x, y)}$-integrable 
for \textit{a.e.} $(x, y)\in\Om\times\square^n$, $\overline{f}$ is in $L^1(\Om\times\square^n)$ and
$f(\cdot, w_h(\cdot))\rightsquigarrow\overline{f}$. 
\end{enumerate}
\end{teo}
\begin{proof}
Assertion (iii) is a straight consequence of Theorems \ref{fund young}(ii) and \ref{multi yuong}.
The integrability properties in assertion (iv) follow by Theorem \ref{fund young}(iii) and 
Remark \ref{oss princ}(v), by noting that $\widehat{\mu}=\widehat{\nu}$.
In order to prove the multiscale convergence, fixed $\varphi\in\ci$ and $\phi\in\cip$,
we define the function $g:\Om\times\square^n\times\R^m\rightarrow\R$ by
\begin{equation*}
g(x, y, z):=\varphi(x)\phi(y)f(x, y, z)\,.
\end{equation*}
The function $g$ is a Caratheodory integrand on $\Om\times T^n\times\R^m$ and $g(\cdot, w_h(\cdot))$ is 
equi-integrable, thus, by Theorems \ref{fund young}(iii) and \ref{multi yuong},
\begin{equation*}
\int_{\Om}g\bigl(x,w_h(x)\bigr)\rightarrow
\int_{\Om}\left(\int_{\square^n\times\R^m}g(x, y, z)d\mu_x(y,z)\right)dx=
\int_{\Om\times\square^n}\varphi(x)\phi(y)\overline{f}(x, y)\,dx\,dy\,.
\end{equation*}
\noindent
Assertion (ii) follows by applying (iv) with $f(x, y, z)=z_\jm$, $\jm=1,\ldots, m$\,.
Finally, by Jensen's inequality and (iii) with $f(x, y, z)=\snorm{z}^q$, we obtain
assertion (i):
\begin{equation*}
\int_{\Om\times\square^n}\snorm{\overline{\nu}(x, y)}^q\,dx\,dy\leq
\int_{\Om\times\square^n}\left(\int_{\R^m}\snorm{z}^q\,d\nu_{(x,y)}(z)\right)dx\,dy
\leq\liminf_{h\rightarrow+\infty}\int_{\Om}\snorm{u_h(x)}^q\,dx<+\infty\,.
\end{equation*}
\end{proof}
\begin{rem}
Actually assertion (ii) is a compactness result about the multiscale convergence:
for every sequence $u_h$ equi-integrable in $L^1$ or bounded in $L^p$, there exists a subsequence 
$u_{h_i}$ which generates a multiscale Young measure and therefore multiscale convergent.
Remember that a bounded sequence in $L^p$ is equi-integrable by H\"older's inequality.
\end{rem}

\vspace{10pt}
We conclude with a basic result about bounded sequences in $W^{1,p}(\Om)$.
\begin{teo}\label{multi}
Given a sequence $u_h$ weakly convergent to $u$ in $W^{1,p}(\Om)$, we have that $u_h\rightsquigarrow u$ and
\begin{equation*}
\g u_h\rightsquigarrow\g u+\sum_{k=1}^n\g_{y^k}u_k
\end{equation*}
for $n$ suitable functions \,$u_k(x, \yk)\in L^p\bigl(\Omega\times\square^{k-1},W_{per}^{1,p}(\square)\bigr)$\,.
\end{teo}
\noindent
The proof can be found in \cite[Theorem 2.6]{All96} and in \cite[Theorem 1.6]{Bald99}.
In the first reference, the idea is to work on the image of $W^{1,2}(\Om)$ under the gradient mapping, 
by characterizing it as the space orthogonal to all divergence-free functions.
Instead in the second reference it is used its characterization as the space of all rotation-free fields.
This last method is simpler and works for general $p$, even if 
only the case $p=2$ is examined in the original statement.
Another proof can be found in \cite{Ped05}.

\section{Continuity results}
\noindent
As first result of this section, we show that it is possible to use in the multiscale convergence a more
complete system of ``test functions", not merely $\psi(x,y)=\varphi(x)\phi(y)$ with $\varphi\in\ci$ and $\phi\in\cip$.
Following Valadier \cite{Val97}, we introduce opportune classes of functions.

\vspace{3pt}
\begin{defi}
A function $\psi:\Om\times\square^n\rightarrow\R$ is said to be \emph{admissible}
if there exist a family $\lbrace X_\delta\rbrace_{\delta>0}$ of compact subsets of $\Om$ 
and a family $\lbrace Y_\delta\rbrace_{\delta>0}$ of compact subsets of $\square$ such that 
$\snorm{\Om\backslash X_\delta}\leq\delta$, $\snorm{\square\backslash Y_\delta}\leq\delta$ 
and $\psi\vert _{X_\delta\times Y_\delta^n}$ is continuous for every $\delta>0$. 
\end{defi}

\begin{rem}
It is not restrictive to suppose that the families $\lbrace X_\delta\rbrace_{\delta>0}$
and $\lbrace Y_\delta\rbrace_{\delta>0}$ are decreasing, \textit{i.e.}, $\delta'\leq\delta$ 
implies $X_{\delta}\subseteq X_{\delta'}$ and $Y_{\delta}\subseteq Y_{\delta'}$. 
Otherwise, it is sufficient to consider the new families $\lbrace\widetilde{X}_\delta\rbrace_{\delta>0}$
and $\lbrace\widetilde{Y}_\delta\rbrace_{\delta>0}$, where
\begin{equation*}
\widetilde{X}_{\delta}:=\bigcap_{i\geq i_{\delta}}X_{2^{-i}}\,, 
\quad \widetilde{Y}_{\delta}:=\bigcap_{i\geq i_{\delta}}Y_{2^{-i}}
\end{equation*}
and $i_{\delta}$ is the minimum positive integer such that $2^{1-i_{\delta}}\leq\delta$.
\end{rem}

\vspace{3pt}
Admissible functions have good measurability properties, as stated in the following lemma.
We omit the easy proof.
\begin{lem}\label{misurabilita}
If $\psi:\Om\times\square^n\rightarrow\R$ is an admissible function, then
there exist a set $X\subseteq\Om$ with $\snorm{\Om\backslash X}=0$ 
and a set $Y\subseteq\square$ with $\snorm{\square\backslash Y}=0$, such that
$\psi\vert _{X\times Y^n}$ is borelian. 
In particular, for every fixed $\e$, the function 
$x\rightarrow \psi\left(x,\langle\frac{x}{\rho_1(\e)}\rangle, \ldots, \langle\frac{x}{\rho_n(\e)}\rangle\right)$ 
is measurable.
\end{lem}

\vspace{3pt}
\begin{defi}
An admissible function $\psi$ is said to be $q$-\emph{admissible}, and we write $\psi\in\mathcal{A}dm^q$,
if there exists a positive function $\alpha\in L^q(\Om)$ such that 
\begin{equation*}
\snorm{\psi(x, y)}\leq \alpha(x) \quad\forall(x, y)\in\Om\times\square^n\,.
\end{equation*}
\end{defi}

\vspace{10pt}
The next theorem proves that it is possible to use $\mathcal{A}dm^q$ as system of test functions.
The proof is very close to \cite[Proposition 5]{Val97}. 
Before we state the following lemma, that can be derived by \cite[Lemma~3.1]{Don83bis}
(see also~\cite[Remark 2.13]{All96}). We use the same definition of $v_h$ given in \eqref{v}.
\begin{lem}\label{disp}
Let $A_k$ be a measurable subset of \;$\square$ for $k=1,\ldots,n$ and let $A:=\prod_{k=1}^nA_k$.
Denoted with $\chi_A$ the characteristic function of $A$, the sequence
$\chi_A(v_h(\cdot))$ converges weakly* to $\snorm{A}$ in $L^{\infty}(\Om)$.
\end{lem}
\begin{teo}\label{continuita1}
Let $u_h$ be a bounded sequence in $L^q(\Om)$, $q\in(1,+\infty]$, generating a multiscale 
Young measure $\nu$ and let $\psi\in\mathcal{A}dm^{q'}$. Then
\begin{equation*}
\lim_{h\rightarrow+\infty}\int_{\Om}\psi\bigl(x, v_h(x)\bigr)u_h(x)dx\\
=\int_{\Om\times\square^n}\psi\left(x, y\right)\overline{\nu}\left(x, y\right)dx\,dy\,.
\end{equation*}
In particular, taking $u_h\equiv1$, we obtain
\begin{equation}\label{caduta}
\lim_{h\rightarrow+\infty}\int_{\Om}\psi\bigl(x, v_h(x)\bigr)dx
=\int_{\Om\times\square^n}\psi\left(x, y\right)dx\,dy\,.
\end{equation}
\end{teo}
\begin{proof}
Let $\delta>0$; by Lusin theorem applied to $\alpha$ and
by definition of admissible function, there exist two compact sets $X\subseteq\Om$ and $Y\subseteq\square$ such that 
$\snorm{\Om\backslash X}\leq\delta$, $\snorm{\square\backslash Y}\leq\delta$ and 
$\psi\vert _{X\times Y^n}$, $\alpha\vert _X$ are continuous. Let $M:=\max_X\alpha$;
by Tietze-Urysohn's theorem, $\psi\vert_{X\times Y^n}$ can be extended 
to a continuous function $\psi_0$ on $\Om\times T^n$ 
with $\snorm{\psi_0(x,y)}\leq M$  for every $(x,y)\in\Om\times\square^n$.
We define on $\Om\times\square^n\times\R$ the functions
\begin{equation*}
f(x, y, z):=\psi(x, y)z \quad\text{and} \quad f_0(x, y, z):=\psi_0(x, y)z\,.
\end{equation*}

\vspace{10pt}
\noindent
With the same definition of $w_h$ given in \eqref{w}, the sequence $f_0(\cdot, w_h(\cdot))$ is equi-integrable because 
$\snorm{f_0(x, w_h(x))}\leq M\snorm{u_h(x)}$. 
By Theorem \ref{fund multi young}(iv), $f_0(\cdot, w_h(\cdot))\rightsquigarrow\overline{f_0}$ and therefore,
by Proposition~\ref{conv debole}, 
$f_0(\cdot, w_h(\cdot))\rightharpoonup\int_{\square^n}\psi_0(\cdot,y)\overline{\nu}(\cdot, y)\,dy$
 weakly in $L^1(\Om)$\,.
This is sufficient to assert that
\begin{equation*}
\lim_{h\rightarrow+\infty}\int_Xf_0\bigl(x, w_h(x)\bigr)dx
=\int_{X\times\square^n}\psi_0\left(x, y\right)\overline{\nu}\left(x, y\right)dx\,dy\,.
\end{equation*}

Now 
\begin{equation*}\begin{split}
\left\vert\int_{\Om\times\square^n}\psi\left(x, y\right)\overline{\nu}\left(x, y\right)dx\,dy\right.
-&\left.\int_{\Om}f\bigl(x,w_h(x)\bigr)\,dx\right\vert
\leq
\snorm{\int_{(\Om\setminus X)\times\square^n}\psi\left(x, y\right)\overline{\nu}\left(x, y\right)dx\,dy}\\
+&
\snorm{\int_{X\times\square^n}\Bigl[\psi\left(x, y\right)-\psi_0\left(x, y\right)\Bigr]\overline{\nu}\left(x, y\right)dx\,dy}\\
+& 
\snorm{\int_{X\times\square^n}\psi_0\left(x, y\right)\overline{\nu}\left(x, y\right)dx\,dy-
       \int_{X}f_0\bigl(x,w_h(x)\bigr)\,dx}\\
+&
\snorm{\int_{X}\Bigl[f_0\bigl(x,w_h(x)\bigr)-f\bigl(x,w_h(x)\bigr)\Bigr]\,dx}
+    
\snorm{\int_{\Om\setminus X}f\bigl(x,w_h(x)\bigr)\,dx}\\
=&\text{I+II+III+IV+V}\,.
\end{split}\end{equation*}

\noindent
We have to show that I, II, IV and V can be made arbitrarily small.
Observe that the function
\begin{equation*}
\gamma_{\text{I}}(x,y):=\alpha(x)\int_{\R}\snorm{z}\,d\nu_{(x,y)}(z)
\end{equation*}
is in $L^1(\Om\times\square^n)$ as consequence of Theorem \ref{fund multi young}(iii), H\"older's inequality
and the $L^q$-boundness of the sequence $u_h$:
\begin{equation*}\begin{split}
\int_{\Om\times\square^n}\gamma_{\text{I}}(x,y)\,dx\,dy\leq&
\liminf_{h\rightarrow+\infty}\int_\Om\alpha(x)\snorm{u_h(x)}\,dx\\
\leq&\;\norm{\alpha}_{L^{q'}(\Om)}\sup_h\norm{u_h}_{L^{q}(\Om)}<+\infty\,.
\end{split}\end{equation*}

\vspace{5pt}
\noindent
The same for $\gamma_{\text{II}}(x,y):=\int_\Om\snorm{z}\,d\nu_{(x,y)}(z)$.
By the absolute continuity of the integral and by the estimates
\begin{align*}
\text{I}\leq&\int_{(\Om\backslash X)\times\square^n}\gamma_{\text{I}}(x,y)\,dx\,dy\\[-15pt]
\intertext{and}
\text{II}=&\int_{X\times(\square^n\backslash Y^n)}\Bigl[\psi(x,y)-\psi_0(x,y)\Bigr]\overline{\nu}\left(x, y\right)dx\,dy\\
\leq&\int_{X\times(\square^n\backslash Y^n)}\gamma_{\text{I}}(x,y)\,dx\,dy
    +M\int_{X\times(\square^n\backslash Y^n)}\gamma_{\text{II}}(x,y)\,dx\,dy\,,
\end{align*}
we obtain that I and II tend to $0$ for $\delta\rightarrow0$.
By using again H\"older's inequality and the $L^q$-boundness of $u_h$, we get
for a suitable positive constant $c$
\begin{align*}
\text{IV}
&\leq\snorm{\int_X\chi_{\square^n\backslash Y^n}(v_h(x))\Bigl[f_0\bigl(x,w_h(x))-f(x,w_h(x)\bigr)\Bigr]dx}\\
&\leq\int_X\chi_{\square^n\backslash Y^n}(v_h(x))\left[\alpha(x)+M\right]\snorm{u_h(x)}dx\\
&\leq c\left(\int_X\chi_{\square^n\backslash Y^n}(v_h(x))
\left[\alpha(x)+M\right]^{q'}dx\right)^{\frac{1}{q'}}\\
\intertext{and}\\[-30pt]
\text{V}
&\leq c\left(\int_{\Om\backslash X}\left[\alpha(x)\right]^{q'}dx\right)^{\frac{1}{q'}}.
\end{align*}

\vspace{10pt}
\noindent
By Lemma \ref{disp}, it follows that
\,$\chi_{\square^n\backslash Y^n}(v_{h}(\cdot))=1-\chi_{Y^n}(v_{h}(\cdot))$  
converges weakly* to $\snorm{\square^n\backslash Y^n}$ and therefore
\begin{equation*}
\int_X\chi_{\square^n\backslash Y^n}(v_h(x))\left[\alpha(x)+M\right]^{q'}dx
\xrightarrow{h\rightarrow\infty}  
\snorm{\square^n\backslash Y^n}\int_X\left[\alpha(x)+M\right]^{q'}dx\,.
\end{equation*}

\vspace{10pt}
\noindent
Hence we conclude that IV and V tend to $0$ for $h\rightarrow\infty$ and $\delta\rightarrow0$\,.
\end{proof}

\vspace{8pt}
\begin{rem}
Let $\psi=\psi(x,\y)$ be a real function on $\Om\times\square^n$ either
continuous in $(y^1,\ldots,y^n)$ and measurable in $x$
or continuous in $(x,y^1,\ldots,y^{k-1},y^{k+1},\ldots,y^n)$
and measurable in $y^k$.
By Scorza-Dragoni theorem (see \cite{Eke99}), $\psi$ is an admissible function.
This is no longer true if one removes the continuity assumption on two variables.
More generally, the invocation of~\eqref{caduta} may be invalid,
as shown in the next two examples. 
The first covers the case $\psi=\psi(x,y)$ ($n=1$)
while the second covers the case $\psi=\psi(y^1, y^2)$.
We remark that in both examples $\psi$ is a Borel function.
See also the example in \cite[Proposition 5.8]{All92}.

\begin{ex}\label{brutto}
Define the Borel sets $A_i:=\bigcup_{j=0}^{i-1}\lbrace(x,y)\in[0,1)^2: \,y=i\,x-j\rbrace$ and $A:=\bigcup_{i=1}^{\infty}A_i$.
Now, in the simple case $d=n=1$, $\rho_1(\e_h)=h^{-1}$ and $\Om=(0,1)$, consider the function $\psi(x,y):=\chi_A(x,y)$.
We have 
\begin{equation*}
\int_0^1\psi(x,\langle h\,x\rangle)\,dx\equiv 1 \quad\text{but} \quad\int_0^1\int_0^1\psi(x,y)\,dx\,dy=0\,.
\end{equation*}
\end{ex}

\begin{ex}
In the case $d=1$, $n=2$, $\rho_1(\e_h)=h^{-1}$, $\rho_2(\e_h)=h^{-2}$ and $\Om=(0,1)$, 
consider the function $\psi(y^1,y^2):=\chi_A(y^1,y^2)$, where $A$ is defined as in the former example.
We have 
\begin{equation*}
\int_0^1\psi(\langle h\,x\rangle,\langle h^2\,x\rangle)\,dx\equiv 1 \quad\text{but} 
\quad\int_0^1\int_0^1\int_0^1\psi(y^1,y^2)\,dx\,dy^1\,dy^2=0\,.
\end{equation*}

\vspace{5pt}
\noindent
Notice that the result of weak* convergence in $L^\infty$ stated in Lemma~\ref{disp}
is not applicable to $A$.
\end{ex}

\end{rem}                                                                                             
\vspace{10pt}
So far we have considered Caratheodory functions $f$ on $\Om\times T^n\times\R^m$. 
As we explained in the introduction, one would like to have a minimal regularity in $(x,\y)$.
For this reason, we introduce an opportune class of integrands and extend to this
Theorem \ref{fund multi young}(iii).

\begin{defi}\label{new}
A function $f:\Om\times\square^n\times\R^m\rightarrow[0,+\infty)$ 
is said to be an \emph{admissible integrand}
if for every $\delta>0$ there exist a compact set $X\subseteq\Om$ with $\snorm{\Om\backslash X}\leq\delta$ and 
a compact set $Y\subseteq\square$ with $\snorm{\square\backslash Y}\leq\delta$, 
such that $f\vert _{X\times Y^n\times\R^m}$ is continuous.
\end{defi}

As in the analogous case for admissible functions (Lemma \ref{misurabilita}),
it is easy to verify the following measurability properties
of admissible integrands.
\begin{lem}
If $f:\Om\times\square^n\times\R^m\rightarrow[0,+\infty)$ is an admissible integrand, then
there exist a set $X\subseteq\Om$ with $\snorm{\Om\backslash X}=0$ 
and a set $Y\subseteq\square$ with $\snorm{\square\backslash Y}=0$, such that
$f\vert _{X\times Y^n\times\R^m}$ is borelian.
In particular, for every fixed $\e$, the function 
$(x,z)\rightarrow f\left(x,\langle\frac{x}{\rho_1(\e)}\rangle, \ldots, \langle\frac{x}{\rho_n(\e)}\rangle,z\right)$ 
is $\f(\Om)\otimes\b(\R^m)$-measurable.
\end{lem}

\vspace{5pt}
\begin{ex}\label{int ammiss}
Let $f:\Om\times\square\times\R^m\rightarrow[0,+\infty)$ be a function such that
\begin{itemize}
\item[(i)]$f(\cdot, y, \cdot)$ is continuous for all $y\in\square$;
\item[(ii)]$f(x, \cdot, z)$ is measurable for all $x\in\Om$ and $z\in\R^m$.
\end{itemize}
By Scorza-Dragoni theorem, $f$ is an admissible integrand. 
Clearly it is possible to replace conditions (i) and (ii) with
\begin{itemize}
\item[(i)']$f(x, \cdot, \cdot)$ is continuous for all $x\in\Om$;
\item[(ii)']$f(\cdot, y,z)$ is measurable for all $y\in\square$ and $z\in\R^m$.
\end{itemize}
\end{ex}

\begin{ex}\label{int ammiss due}  
Let $f:\Om\times\square^n\times\R^m\rightarrow[0,+\infty)$ be a function of the type
\begin{equation*}\begin{split}
f(x,\y,z)=\prod_{k=1}^n \chi_{P_k}(y^k)\,f_1(x,\y,z)&\\+\Bigl[1-\prod_{k=1}^n \chi_{P_k}(y^k)\Bigr]&\,f_2(x,\y,z)\,,
\end{split}\end{equation*}
where $P_k$ ($k=1,\ldots,n$) is a  measurable subset of $\square$ and
$f_j$ ($j=1,2$) is a non-negative function on $\Om\times\square^n\times\R^m$ such that
\begin{itemize}
\item[(i)]$f_j$ is continuous in $(y^1,\ldots,y^n,z)$;
\item[(ii)]$f_j$ is measurable in $x$.
\end{itemize}
By Scorza-Dragoni theorem for every $\delta>0$ there exists a compact set $X\subseteq\Om$
such that the functions $f_1$ and $f_2$ are continuous on $X\times\square^n\times\R^m$. 
By applying Lusin theorem to each $\chi_{P_k}$, we obtain that $f$ is an admissible integrand.
Obviously, the conditions (i) and (ii) can be replaced by
\begin{itemize}
\item[(i)']$f_j$ is continuous in $(x,y^1,\ldots,y^{k-1},y^{k+1},\ldots,y^n,z)$;
\item[(ii)']$f_j$ is measurable in $y^k$.
\end{itemize}
\end{ex}

\vspace{6pt}
\begin{teo}\label{continuita2}
Let $u_h$ be a bounded sequence in $L^q(\Om,\R^m)$, $q\in[1,+\infty)$, generating a multiscale 
Young measure $\nu$ and let $f:\Om\times\square^n\times\R^m\rightarrow[0,+\infty)$  be an 
admissible integrand satisfying the $q$-growth condition
\begin{equation*}
f\left(x, \y, z\right)\leq c\left(1+\snorm{z}^q\right)
\end{equation*}
for some $c>0$ and for all $\left(x, \y, z\right)\in\Om\times\square^n\times\R^m$\,. Then
\begin{equation}\label{master}
\liminf_{h\rightarrow+\infty}\int_{\Om}
f\bigl(x, w_h(x)\bigr)dx\geq\int_{\Om\times\square^n}\overline{f}(x, \y)\,dx\,dy^1\ldots dy^n ,
\end{equation}
where as usual
\begin{equation*}
\overline{f}(x, \y):=\int_{\R^m}f(x, \y, z)\;d\nu_{(x, \y)}(z)\,.
\end{equation*}
\end{teo}
\begin{proof}
Assume initially that, in addition, 
\begin{equation}\label{limitatezza}
f(x, y, z)=0 \quad\text{if} \; \snorm{z}\geq r
\end{equation}
for a fixed $r>0$. By the admissibility condition, for every
$\delta>0$ there exist a compact set $X\subseteq\Om$ and  
a compact set $Y\subseteq\square$ such that $\snorm{\Om\backslash X}\leq\delta$, 
$\snorm{\square\backslash Y}\leq\delta$ and $f\vert _{X\times Y^n\times\R^m}$ is continuous.
By Tietze-Urysohn's theorem, $f\vert_{X\times Y^n\times\R^m}$ can be extended 
to a continuous function $f_0$ on $\Om\times T^n\times\R^m$ 
with $0\leq f_0(x, y, z)\leq M$ for every $(x, y, z)\in\Om\times\square^n\times\R^m$,
where $M:=\max f$ on $X\times Y^n\times\R^m$.
Notice that, by the $q$-growth condition, $M\leq c\left(1+r^q\right)$.
Obviously $f_0\left(\cdot,w_h(\cdot)\right)$ is equi-integrable and so, by Theorem \ref{fund multi young}(iv) and
by Proposition \ref{conv debole},
\begin{equation*}
\lim_{h\rightarrow+\infty}\int_{\Om} f_0\bigl(x,w_h(x)\bigr)\,dx=
\int_{\Om\times\square^n}\overline{f_0}(x, y)\,dx\,dy\,.
\end{equation*}

\vspace{8pt}
\noindent
For a suitable subsequence $h_i$ 
\begin{equation*}
\lim_{i\rightarrow+\infty}\int_{\Om}f\bigl(x,w_{h_i}(x)\bigr)dx=
\liminf_{h\rightarrow+\infty}\int_{\Om}f\bigl(x,w_h(x)\bigr)dx.
\end{equation*} 
We can write
\begin{equation*}\begin{split}
\lim_{i\rightarrow+\infty}\int_{\Om}f\bigl(x,w_{h_i}(x)\bigr)\,dx\,-&\int_{\Om\times\square^n}\overline{f}(x, y)\,dx\,dy
=\lim_{i\rightarrow+\infty}\int_{\Om}\Bigl[f\bigl(x,w_{h_i}(x)\bigr)-f_0\bigl(x,w_{h_i}(x)\bigr)\Bigr]dx\\
&\,+\Biggr[\lim_{i\rightarrow+\infty}\int_{\Om}f_0\bigl(x,w_{h_i}(x)\bigr)\,dx-\int_{\Om\times\square^n}\overline{f_0}(x, y)\,dy\,dx\Biggl]\\
&+\int_{\Om\times\square^n}\Bigl[\overline{f_0}(x, y) -\overline{f}(x, y)\Bigr]\,dx\,dy
=\text{I+II+III}\,.
\end{split}\end{equation*}
Let us check that the negative part of \,I and III can be made arbitrarily small. 

Firstly, by Lemma \ref{disp},
\begin{equation*}\begin{split}
\text{I}=&
\lim_{i\rightarrow+\infty}\int_{\Om\backslash X}\Bigl[f\bigl(x,w_{h_i}(x)\bigr)-f_0\bigl(x,w_{h_i}(x)\bigr)\Bigr]dx\\
  &+\lim_{i\rightarrow+\infty}\int_{X}
  \chi_{\square^n\backslash Y^n}(v_{h_i}(x))\Bigl[f\bigl(x,w_{h_i}(x)\bigr)-f_0\bigl(x,w_{h_i}(x)\bigr)\Bigr]dx\\
  \geq&-M\snorm{\Om\backslash X}\,-\lim_{i\rightarrow+\infty}M\int_{X}\chi_{\square^n\backslash Y^n}(v_{h_i}(x))\,dx\\
  \geq&-M\bigl(\snorm{\Om\backslash X}+\snorm{X}\snorm{\square^n\backslash Y^n}\bigr)
  \geq-c\left(1+r^q\right)\bigl(\delta+n\snorm{\Om}\delta\bigr)\,.
\end{split}\end{equation*}

\vspace{8pt}
\noindent
Now, observe that the function
\begin{equation*}
\gamma(x,y):=\int_{\R^m}\Bigl[f(x, y, z)-f_0(x, y, z)\Bigr]\,d\nu_{(x,y)}(z)
\end{equation*}
is in $L^1(\Om\times\square^n)$ as consequence of Theorem \ref{fund multi young}(iii):
\begin{equation*}\begin{split}
\int_{\Om\times\square^n}\snorm{\gamma(x,y)}\,dx\,dy\leq&
\int_{\Om\times\square^n}\left[M+c+c\int_{\R^m}\snorm{z}^q d\nu_{(x,y)}(z)\right]\,dx\,dy\\
\leq&\;(M+c)\snorm{\Om}+c\,\liminf_{h\rightarrow+\infty}\int_{\Om}\snorm{u_h(x)}^q\,dx
<+\infty\,.
\end{split}\end{equation*}
By the absolute continuity of the integral and by the equality
\begin{equation*}
\text{III}=\int_{(\Om\times\square^n)\setminus(X\times Y^n)}\gamma(x,y)\,dx\,dy\,,
\end{equation*}
we obtain that III tends to $0$ for $\delta\rightarrow0$.
This concludes the first part of the proof.\\[10pt]
In order to remove assumption \eqref{limitatezza} we consider,
for $k\in\NP$, the functions \,$p_k\in C_0(\R^m)$ defined by
\begin{equation*}
p_k\left(z\right):=
\begin{cases}
1 & \textrm{if} \quad \snorm{z}\leq k\\
1+k-\snorm{z} & \textrm{if} \quad k\leq\snorm{z}\leq k+1\\ 
0 & \textrm{if} \quad  \snorm{z}\geq k+1
\end{cases}
\end{equation*}

\vspace{10pt}
\noindent
and the functions $f_k(x, y, z):=p_k(z)f(x, y, z)$\,. By applying the first part of the theorem, we have

\begin{equation*}
\liminf_{h\rightarrow+\infty}\int_{\Om}f\bigl(x,w_h(x)\bigr)dx\geq
\liminf_{h\rightarrow+\infty}\int_{\Om}f_k\bigl(x,w_h(x)\bigr)dx\geq
\int_{\Om\times\square^n}\overline{f_k}(x, y)\,dx\,dy\,.
\end{equation*}

\vspace{8pt}
\noindent
By noting that $f_k$ is increasing and that $f_k(x, y, \cdot)\rightarrow f(x, y, \cdot)$ \textit{a.e.} in $\R^m$
for every fixed $(x,y)\in\Om\times\square^n$, we deduce from the monotone convergence theorem that 
$\overline{f_k}\rightarrow\overline{f}$ \textit{a.e.} in $\Om\times\square^n$. 
The sequence $\overline{f_k}$ is increasing so, again from monotone convergence theorem, 
\begin{equation*}
\int_{\Om\times\square^n}\overline{f_k}(x, y)\,dx\,dy\xrightarrow{k\rightarrow\infty} 
\int_{\Om\times\square^n}\overline{f}(x, y)\,dx\,dy\,.
\end{equation*}
\end{proof}

\begin{rem}
Lower semicontinuity property \eqref{master} is not true if $f$ is only borelian.
For instance, consider the function 
$f(x,y,z):=[1-\psi(x,y)]\snorm{z}^p$, where $\psi$ is defined as in Example~\ref{brutto}.
\end{rem}

\section{Gamma-convergence}
\noindent
In the present section we examine the multiperiodic homogenization of nonlinear convex functionals
by means of the $\Gamma$-convergence combined with the multiscale Young measures.

Before we recall the definition of $\Gamma$-convergence,
referring to \cite{Br98} and \cite{Dal93} for an exposition of the main properties.

\begin{defi}
Let $(U,\tau)$ be a topological space satisfying the first countability axiom and
$F_h$, $F$ functionals from $U$ to $[-\infty, +\infty]$;
we say that $F$ is the $\Gamma(\tau)$-limit of the sequence $F_h$ or that $F_h$ $\Gamma(\tau)$-converges to $F$,
and write
\begin{equation*}
F=\Gamma(\tau)\hbox{-}\negthinspace\lim_{h\rightarrow+\infty}F_h,
\end{equation*}
if for every $u\in U$ the following conditions are satisfied:

\begin{equation}\label{disug inf}
\hspace*{-5pt}F(u)\leq\inf\Bigl\lbrace\liminf_{h\rightarrow+\infty}F_h(u_h) : u_h\xrightarrow{\tau}u\Bigr\rbrace
\end{equation}
and
\begin{equation}\label{disug sup}
F(u)\geq\inf\Bigl\lbrace\limsup_{h\rightarrow+\infty}F_h(u_h) : u_h\xrightarrow{\tau}u\Bigr\rbrace\,.
\end{equation}
\end{defi}

\vspace{8pt}
We can extend the definition of $\Gamma$-convergence to families depending on a parameter $\e>0$.
\begin{defi}
For every $\e>0$, let ${F_\e}$ be a functional from $U$ to $[-\infty, +\infty]$.
We say that $F$ is the $\Gamma(\tau)$-limit of the family $F_\e$, and write
\begin{equation*}
F=\Gamma(\tau)\hbox{-}\negthinspace\lim_{\e\rightarrow0^+}F_\e,
\end{equation*}
if we have for every sequence $\e_h\rightarrow0^+$
\begin{equation*}
F=\Gamma(\tau)\hbox{-}\negthinspace\lim_{h\rightarrow+\infty}F_{\e_h}.
\end{equation*}
\end{defi}

\vspace{10pt}
Throughout this section, we work in the space $L^p(\Om,\R^m)$ endowed with
the strong topology. As pointed out in the introduction, we consider a non-negative function 
$f=f\left(x,\y, z\right)$ on $\Om\times\square^n\times\M^{m\times d}$
satisfying Assumptions \ref{ass1}, \ref{ass2} and \ref{ass3}. 

\vspace{10pt}
We fully characterize the $\Gamma(L^p)$-limit of the family 
$F_\e:L^p(\Om,\R^m)\rightarrow[0, +\infty]$ where the functionals are defined by
\begin{equation*}
F_\e (u):=
\begin{cases}
\displaystyle\int_\Om f\Bigl(x,\langle\frac{x}{\rho_1(\e)}\rangle, \ldots, \langle\frac{x}{\rho_n(\e)}\rangle, \g u(x)\Bigr)dx
& \textrm{if} \quad u\in W^{1,p}(\Om,\R^m),\\
&\\
+\infty & \textrm{otherwise.}
\end{cases}
\end{equation*}

\vspace{10pt}
\noindent
Precisely, this is our main result.
\begin{teo}\label{main1}
The family $F_\e$ \,$\Gamma(L^p)$-converges and its $\Gamma(L^p)$-limit \,$F_{hom}:L^p(\Om,\R^m)\rightarrow[0, +\infty]$
is given by
\begin{equation*}
F_{hom} (u)=
\begin{cases}
\displaystyle\int_\Om f_{hom}\bigl(x,\g u(x)\bigr)dx & \textrm{if} \quad u\in W^{1,p}(\Om,\R^m),\\
&\\
+\infty & \textrm{otherwise},
\end{cases}
\end{equation*}

\vspace{6pt}
\noindent
where $f_{hom}$ is obtained by the following cell problem
\begin{equation*}
f_{hom}\left(x,z\right):=
\inf_{\phi\in\Phi}\int_{\square^n} f\Bigl(x,y, z+\sum_{k=1}^n\g_{y^k}\phi_k(\yk)\Bigr)dy 
\end{equation*}
with the space $\Phi$ defined by
\begin{equation*}
\Phi:=\prod_{k=1}^n\Phi_k \quad\text{and} \quad\Phi_k:=L^p\bigl(\square^{k-1},W_{per}^{1,p}(\square,\R^m)\bigr)\,.
\end{equation*}
\end{teo}

\vspace{8pt}
\begin{rem}\label{sup=supreg}
(i) Using the $p$-growth condition of $f$ and a density argument, it can be shown that
\begin{equation*}
f_{hom}\left(x,z\right)=
\inf_{\phi\in\Phi_{reg}}\int_{\square^n} f\Bigl(x,y, z+\sum_{k=1}^n\g_{y^k}\phi_k(\yk)\Bigr)\,dy\,,
\end{equation*}

\vspace{-3pt}
\noindent
where
\begin{equation*}
\Phi_{reg}:=\prod_{k=1}^n\Phi_{k, reg} \quad\text{and}
\quad\Phi_{k, reg}:=C^1\Bigl(\overline{\square}^{k-1},C_{per}^1(\square,\R^m)\Bigr)\,.
\end{equation*}

\vspace{3pt}
\noindent
(ii) For every $\delta>0$ there exists a compact set $X\subseteq\Om$ with $\snorm{\Om\backslash X}\leq\delta$ 
such that the restriction of $f$ to $X\times\square^n\times\M^{m\times d}$ 
is continuous in $(x,z)$ for \textit{a.e.} $(\y)\in\square^n$ and so
$f_{hom}$ is lower semicontinuous on $X\times\M^{m\times d}$.
In particular $f_{hom}$ is $\f(\Om)\otimes\b(\M^{m\times d})$-measurable.\\
(iii) The convexity, the $p$-coerciveness and the $p$-growth condition on $f$ give the 
corresponding properties for the function $f_{hom}$.
In particular $F_{hom}$ is continuous on $W^{1,p}(\Om,\R^m)$, endowed with the strong topology.
\end{rem}

\vspace{6pt}
Before proving the theorem, we state a series of lemmas.
Only for simplicity of notations, we restrict ourselves to the case $m=1$.
Fixed a sequence $\e_h\rightarrow0^+$, we use 
for $v_h$ the same definition given in \eqref{v}.
\begin{lem}\label{inf lemma}
Let $u_h$ be a sequence converging weakly in $W^{1,p}(\Om)$ to a function $u$. Then

\begin{equation*}
\liminf_{h\rightarrow+\infty}\int_\Om f\bigl(x, v_h(x), \g u_h(x)\bigr)dx\geq
\int_\Om f_{hom}\bigl(x,\g u(x)\bigr)dx\,.
\end{equation*}
\end{lem}
\begin{proof} 
For a suitable subsequence $h_i$\,,
\begin{equation*}
\lim_{i\rightarrow+\infty}\int_\Om f\bigl(x, v_{h_i}(x), \g u_{h_i}(x)\bigr)dx=
\liminf_{h\rightarrow+\infty}\int_\Om f\bigl(x, v_h(x), \g u_h(x)\bigr)dx\,.
\end{equation*}
Refining the subsequence if necessary, we can suppose that $u_{h_j}$ generates a multiscale
Young measure $\nu$. By Theorem \ref{continuita2} and Jensen's inequality
\begin{align*}
\lim_{i\rightarrow+\infty}\int_\Om f\bigl(x, v_{h_i}(x), \g u_{h_i}(x)\bigr)dx
\geq&\int_{\Om\times\square^n}\int_{\R^d}f(x, y, z)\,d\nu_{(x,y)}\,dx\,dy\\
\geq&\int_{\Om\times\square^n}f\Bigl(x, y, \int_{\R^d}z\,d\nu_{(x,y)}\Bigr)\,dx\,dy\\[-10pt]
\intertext{and by Theorems \ref{fund multi young}(ii) and \ref{multi}}\\[-25pt]
\geq&\int_{\Om\times\square^n}f\Bigl(x, y, \g u(x)+\sum_{k=1}^n\g_{y^k}u_k(x, \yk)\Bigr)\,dx\,dy\\
\geq&\int_\Om f_{hom}\bigl(x, \g u(x)\bigr)dx\,.
\end{align*}
\end{proof}

\begin{lem}\label{lemma convessita}
Let $f:\R^d\rightarrow\R$ be a convex function, such that for every $z\in\R^d$
\vspace{3pt}
\begin{equation}\label{stima crescita}
\snorm{f(z)}\leq c\left(b+\snorm{z}\right)^p ,
\end{equation}
\noindent
where b and c are positive constants. Then, for all $z_1, z_2\in\R^d$
\vspace{3pt}
\begin{equation}\label{stima convessa}
\snorm{f(z_1)-f(z_2)}\leq c\,d\,(1+2^p)\bigl(b+\snorm{z_1}+\snorm{z_2}\bigr)^{p-1}\snorm{z_1-z_2}\,.
\end{equation}
\end{lem}
\noindent
The proof can be derived from \cite[Lemma 5.2]{Giu03}.  We observe that in \eqref{stima convessa}
the estimate depends only by the costants \textit{b, c} of growth condition \eqref{stima crescita}
and not by the particular function $f$.

\vspace{8pt}
\begin{lem}\label{sup lemma}
Let $u\in W^{1,p}(\Om)\cap C^1(\Om)$. Then
\begin{equation}\label{ineq}
\inf_{u_h\rightarrow u}\Bigl\lbrace\limsup_{h\rightarrow+\infty}F_{\e_h}(u_h)\Bigr\rbrace
\leq\inf_{\psi\in\Psi}\int_{\Om\times\square^n}f\Bigl(x, y, \g u(x)+\sum_{k=1}^n\g_{y^k}\psi_k(x, \yk)\Bigr)\,dx\,dy\,,
\end{equation}
where the $\inf$'s are made respectively on the sequences $u_h$ that converge strongly in $L^p(\Om)$ to $u$
and on the space $\Psi$ defined by
\begin{equation*}
\Psi:=\prod_{k=1}^n\Psi_k \quad\text{and}
\quad\Psi_k:=C^1\left(\overline{\Om}\times\overline{\square}^{k-1},C_{per}^1(\square)\right)\,.
\end{equation*}
\end{lem}
\begin{proof}
Given an arbitrary function $\psi=(\psi_1,\ldots ,\psi_n)\in\Psi$, consider the sequence 
\begin{equation*}
u_h(x):=u(x)+\sum_{k=1}^n \rho_k(\e_h)\psi_k\bigl(x,v^k_h(x)\bigr)\,,
\end{equation*}
where we used the short notation 
$v^k_h(x):=\left(\langle\frac{x}{\rho_1(\e_h)}\rangle, \ldots, \langle\frac{x}{\rho_k(\e_h)}\rangle\right)\,.$
We have $u_h\rightarrow u$ strongly in $L^p(\Om)$ and $\g u_h=\g u+\sum_{k=1}^n\g_{y^k}\psi_k+r_h$, with
$r_h\rightarrow 0$ strongly in $L^p(\Om,\R^d)$.

\vspace{5pt}
The function $g:\Om\times\square^n\rightarrow\R$ defined by
\begin{equation*}
g(x,y):=f\Bigl(x, y, \g u(x)+\sum_{k=1}^n\g_{y^k}\psi_k(x, \yk)\Bigr)
\end{equation*}
is admissible. Actually $g\in\mathcal{A}dm^1$, as evident by the estimate obtained  through the 
$p$-growth condition:

\begin{equation*}
\snorm{g(x,y)}\leq c_2\Bigl[1+(n+1)^{p-1}
\Bigl(\snorm{\g u(x)}^p+\sum_k M_k^p\Bigr)\Bigr]\,,
\end{equation*}
where $M_k:=\sup_{\Om\times\square^k}\vert\g_{y^k}\psi_k\vert$. 

\vspace{5pt}
By Lemma \ref{lemma convessita},
the following inequality holds for some positive constants $b$, $c$\,:
\begin{equation*}
\Bigl\vert g\bigl(x,v_h(x)\bigr)-f\bigl(x,v_h(x),\g u_h(x)\bigr)\Bigr\vert
\leq c \snorm{r_h(x)}\Bigl(b+\snorm{\g u(x)}^{p-1}+\snorm{r_h(x)}^{p-1}\Bigr)\,.
\end{equation*}
By integrating over $\Om$, from H\"older's inequality
we obtain, for another positive constant $c'$,
\begin{equation*}
\int_{\Om}\Bigl\vert g\bigl(x,v_h(x)\bigr)-f\bigl(x,v_h(x),\g u_h(x)\bigr)\Bigr\vert\,dx
\leq c' \int_{\Om}\snorm{r_h(x)}^p dx
\end{equation*}
and thus Theorem \ref{continuita1} gives
\begin{equation*}\begin{split}
\lim_{h\rightarrow+\infty}&\int_{\Om}f\bigl(x,v_h(x),\g u_h(x)\bigr)\,dx
=\lim_{h\rightarrow+\infty}\int_{\Om}g\bigl(x, v_h(x)\bigr)dx\\
=&\int_{\Om\times\square^n}g\left(x, y\right)dx\,dy
=\int_{\Om\times\square^n}f\Bigl(x, y, \g u(x)+\sum_{k=1}^n\g_{y^k}\psi_k(x, \yk)\Bigr)\,dx\,dy\,.
\end{split}\end{equation*}
\end{proof}

\vspace{3pt}
\begin{defi}
We say that $\Lambda\subseteq L^1(\Om)$ is an \emph{inf-stable family} if, given 
$\lbrace\lambda_1,\ldots ,\lambda_N\rbrace\subseteq\Lambda$ and
$\lbrace\varphi_1,\ldots ,\varphi_N\rbrace\subseteq C^1\bigl(\overline{\Om},[0,1]\bigr)$,
with $\sum_{\jm=1}^N\varphi_\jm=1$ and $N\in\NP$, there exists a $\lambda\in\Lambda$ such that 
\begin{equation*}
\lambda\leq\sum_{\jm=1}^N\varphi_\jm\lambda_\jm\,.
\end{equation*}
\end{defi}
\begin{lem}\label{lemma commutativita}
Let $\Lambda$ be an inf-stable family of non-negative integrable functions on $\Om$.
If for every $\delta>0$ there exists a compact set $X_{\delta}\subseteq\Om$ such that 
$\snorm{\Om\backslash X_{\delta}}\leq\delta$ 
and $\lambda\vert _{X_{\delta}}$ is continuous for each $\lambda\in\Lambda$,
then the function $\inf_{\lambda\in\Lambda}\lambda$ is measurable and the following commutation property holds:
\begin{equation}\label{inf comm}
\inf_{\lambda\in\Lambda}\int_{\Om}\lambda(x)\,dx = \int_{\Om}\inf_{\lambda\in\Lambda}\lambda(x)\,dx\,.
\end{equation}
\end{lem}
\noindent
This lemma can be derived by \cite[Lemma 4.3]{Bou93} (see also \cite{Haf03}),
by noting that for every $\delta>0$ 
$\inf_{\lambda\in\Lambda}\lambda=\essinf_{\lambda\in\Lambda}\lambda$ on $X_{\delta}$.
Anyway, we prefer to give a simple direct proof.
\begin{proof}
Firstly we observe that for every $\delta>0$ the function
$\inf_{\lambda\in\Lambda}\lambda$ is lower semicontinuous on~$X_{\delta}$.
In particular  $\inf_{\lambda\in\Lambda}\lambda$ is measurable.
By applying the Lindel\"of theorem to each family $\lbrace E_{\delta}^{\lambda}\rbrace_{\lambda\in\Lambda}$,
where
\begin{equation*}
E_{\delta}^{\lambda}:=\bigl\lbrace(x,t)\in X_{\delta}\times\R : \lambda(x)<t\bigr\rbrace,
\end{equation*}

\vspace{5pt}
\noindent
we can find a sequences $\lambda_i$ in $\Lambda$ such that
\begin{equation*}
\inf_{\lambda\in\Lambda}\lambda(x)=\inf_i\lambda_i(x)  \quad \text{for \textit{a.e.}} \;x\in\Om\,.
\end{equation*}
Fixed $N\in\NP$ and $\zeta>0$, we choose a $\delta>0$ such that 
$\sum_{j=1}^N\int_{\Om\backslash X_{\delta}}\lambda_j\leq\zeta$.
By the continuity property of the elements $\lambda\in\Lambda$, the sets 
\begin{equation*}
A_i:=\left\lbrace x\in X_{\delta} : \lambda_i(x)<\inf_{1\leq j \leq N}\lambda_j(x)+\zeta\right\rbrace
\end{equation*}
are open in $X_{\delta}$. Notice that $X_{\delta}=\bigcup_{i=1}^\infty A_i$.
For every $i\in\NP$, let $B_i$ be a open subset of $\R^d$ for which $B_i\cap X_{\delta}=A_i$ and
let $\lbrace\varphi_i\rbrace_i\subseteq C^1\bigl(\overline{\Om},[0,1]\bigr)$ be
a partition of unity subordinate to $\lbrace B_i\rbrace_i$. By the inf-stability property,
there exists a $\lambda\in\Lambda$ such that $\lambda\leq\sum_{j=1}^N\varphi_j\lambda_j$. We have
\begin{equation*}\begin{split}
\int_\Om\lambda(x)\,dx=&\int_{\Om\backslash X_{\delta}}\lambda(x)\,dx + \int_{X_{\delta}}\lambda(x)\,dx\\
\leq&\,\sum_{j=1}^N\int_{\Om\backslash X_{\delta}}\varphi_j(x)\,\lambda_j(x)\,dx
+\sum_{i=1}^\infty\int_{X_{\delta}}\varphi_i(x)\lambda_i(x)\,dx\\
\leq&\,\,\zeta+\int_\Om\inf_{1\leq j\leq N}\lambda_j(x)\,dx + \zeta\snorm{\Om}\,.
\end{split}\end{equation*}
Being $N$ and $\zeta$ arbitrary, the claim follows.
\end{proof}

\vspace{10pt}
\noindent
We are now ready to assemble a proof of Theorem \ref{main1}.

Let $u_h\rightarrow u$ in $L^p(\Om)$. We want to show that $\liminf F_{\e_h}(u_h)\geq F_{hom}(u)$.
In this way inequality~\eqref{disug inf} will be proved.
If $\liminf F_{\e_h}(u_h)=+\infty$, there is nothing to prove, so
we can assume $\liminf F_{\e_h}(u_h)<+\infty$. For a suitable subsequence $h_i$\,,
\begin{equation*}
\lim_{i\rightarrow+\infty}F_{\e_{h_i}}(u_{h_i})=
\liminf_{h\rightarrow+\infty}F_{\e_h}(u_h)\,.
\end{equation*}
For $i$ large enough, $F_{\e_{h_i}}(u_{h_i})$ is finite and therefore, by the definition of $F_\e$,
$u_{h_i}\in W^{1,p}(\Om)$. 
Due to the $p$-coerciveness hypothesis on $f$, we can infer that $\g u_{h_i}$ is bounded in $W^{1,p}(\Om)$.
Refining the subsequence if necessary, we can suppose that $u_{h_j}$ converges weakly in $W^{1,p}(\Om)$ to $u$ and 
thus we can apply Lemma \ref{inf lemma}.

It remains to check inequality \eqref{disug sup}.
If $u\in L^p(\Om)\setminus W^{1,p}(\Om)$, then $F_{hom}(u)=+\infty$ and the inequality
is obvious, while if $u\in W^{1,p}(\Om)$, then we can apply Lemma \ref{sup lemma}
and, as in \cite[Theorem 3.3]{Bou04}, Lemma \ref{lemma commutativita}.
In view of the density of $W^{1,p}(\Om)\cap C^1(\Om)$ in $W^{1,p}(\Om)$ and of the continuity of $F_{hom}$,
by a standard diagonalization argument, it is not restrictive to assume that $u\in W^{1,p}(\Om)\cap C^1(\Om)$.

For every $\psi=(\psi_1,\ldots ,\psi_n)\in\Psi$, define the function
\begin{equation*}
\lambda_{\psi}(x):=\int_{\square^n}f\Bigl(x, y, \g u(x)+\sum_{k=1}^n\g_{y^k}\psi_k(x, \yk)\Bigr)\,dy\,.
\end{equation*}
We claim that the family $\Lambda:=\left\lbrace\lambda_\psi : \psi\in\Psi\right\rbrace$
satisfies the hypotheses of Lemma \ref{lemma commutativita}.
In fact, from the $p$-growth condition on $f$, it is easy to show that each function in $\Lambda$ is 
integrable on $\Om$. Moreover, by Remark \ref{sup=supreg}(ii), for every $\delta>0$
there exists a compact set $X\subseteq\Om$ with $\snorm{\Om\backslash X}\leq\delta$ 
such that $\lambda_{\psi}$ is continuous on $X$ for each $\psi\in\Psi$.
It remains to prove the inf-stability.

Given $\lbrace\psi^{(1)},\ldots ,\psi^{(N)}\rbrace\subseteq\Psi$ and
$\lbrace\varphi_1,\ldots ,\varphi_N\rbrace\subseteq C^1\bigl(\overline{\Om},[0,1]\bigr)$,
with $\sum_{\jm=1}^N\varphi_\jm=1$ and $N\in\NP$, consider the function 
\begin{equation*}
\psi:=\Biggl(\sum_{\jm=1}^N\varphi_\jm\psi^{(\jm)}_1,\ldots ,\sum_{\jm=1}^N\varphi_\jm\psi^{(\jm)}_n\Biggr)\in\Psi\,.
\end{equation*}
Thanks to the convexity of $f$, we have $\lambda_{\psi}\leq\sum_{\jm=1}^N\varphi_\jm\lambda_{\psi^{(\jm)}}$:
\begin{equation*}\begin{split}
\lambda_{\psi}(x)
=&\int_{\square^n}f\Bigl(x, y, \g u(x)+\sum_{\jm=1}^N\sum_{k=1}^n\g_{y^k}\varphi_\jm(x)\psi^{(\jm)}_k(x, \yk)\Bigr)\,dy\\
=&\int_{\square^n}f\Biggl(x, y, \sum_{\jm=1}^N\varphi_\jm(x)\Bigl(\g u(x)+
                             \sum_{k=1}^n\g_{y^k}\psi^{(\jm)}_k(x, \yk)\Bigr)\Biggr)\,dy\\
\leq&\sum_{\jm=1}^N\varphi_\jm(x)\int_{\square^n}f\Bigl(x, y, \g u(x)+\sum_{k=1}^n\g_{y^k}\psi^{(\jm)}_k(x, \yk)\Bigr)\,dy
=\sum_{\jm=1}^N\varphi_\jm(x)\lambda_{\psi^{(\jm)}}(x)\,.
\end{split}\end{equation*}

\vspace{8pt}
\noindent
Finally, by inequality \eqref{ineq}, equality \eqref{inf comm} and Remark \ref{sup=supreg}(i), 
\begin{align*}
\inf_{u_h\rightarrow u}&\Bigl\lbrace\limsup_{h\rightarrow+\infty}F_{\e_h}(u_h)\Bigr\rbrace
\leq\int_{\Om}\inf_{\psi\in\Psi}
\left(\int_{\square^n}f\Bigl(x, y, \g u(x)+\sum_{k=1}^n\g_{y^k}\psi_k(x, \yk)\Bigr)\,dy\right)dx\\
&\leq\int_{\Om}\inf_{\phi\in\Phi_{reg}}
\left(\int_{\square^n}f\Bigl(x, y, \g u(x)+\sum_{k=1}^n\g_{y^k}\phi_k(\yk)\Bigr)\,dy\right)dx
=\int_\Om f_{hom}\bigl(x, \g u(x)\bigr)\,dx\,.   
\end{align*}
The proof is complete.
\begin{flushright}$\blacklozenge$\end{flushright}

\vspace{7pt}
Assumpion \ref{ass3} cannot be weakened too much: even if $f$ is a Borel
function, $\Gamma$-convergence Theorem~\ref{main1} may be not applicable,
as shown in the following examples (see also~\cite[Example~3.1]{Don83}).
\begin{ex}\label{final uno}
Let $A_i$ be the Borel sets defined as in Example \ref{brutto} and 
let $B:=\bigcup_{i=1}^{\infty}A_{2i}$. 
Notice that $\bigcup_{i\neq j}(A_i\cap A_j)$ is countable.
In the case $d=n=m=1$, $\rho_1(\e_h)=h^{-1}$ and $\Om=(0,1)$, 
consider the Borel function $f(x,y,z):=[2-\chi_B(x,y)]\snorm{z}^p$.
We remark that $f$ satisfies only Assumptions \ref{ass1} and \ref{ass2}.
We have for every $u\in W^{1,p}((0,1))$
\begin{equation*}
F_h(u)=\int_0^1 f\bigl(x,\langle h\,x\rangle, \g u(x)\bigr)\,dx=
\begin{cases}
\displaystyle\int_0^1\snorm{\g u(x)}^p\,dx
& \textrm{if} \quad h\equiv 0\negthinspace\mod 2,\\
&\\
\displaystyle2\int_0^1\snorm{\g u(x)}^p\,dx
& \textrm{if} \quad h\equiv 1\negthinspace\mod 2.\\
\end{cases}
\end{equation*}
Clearly the sequence $F_h$ is not $\Gamma$-convergent in $L^p((0,1))$ with respect
to the strong topology.
\end{ex}

\begin{ex}\label{final due}
Let $d=m=1$,  $n=2$, $\rho_1(\e_h)=h^{-1}$, $\rho_2(\e_h)=h^{-2}$ and $\Om=(0,1)$.
Consider the Borel function
$f(x,y^1,y^2,z):=[2-\chi_B(y^1,y^2)]\snorm{z}^p$,
where $B$ is defined as in the former example.
Even if $f$ does not depend by $x$ and satisfies Assumptions \ref{ass1} and \ref{ass2},
the sequence $F_h$ is not $\Gamma(L^p)$-convergent.
\end{ex}

\section{Iterated homogenization}
\noindent
The homogenized function $f_{hom}$ can be obtained also by the following iteration:
\begin{equation*}\begin{split}
&f_{hom}^{[n]}\left(x,y^1, \ldots,y^{n-1}, z\right):=
\inf_{\phi\in W^{1,p}_{per}(\square,\R^m)}
\int_\square f\left(x,y^1, \ldots,y^n, z+\g \phi(y^n)\right)dy^n,\\
&f_{hom}^{[n-1]}\left(x,y^1, \ldots,y^{n-2}, z\right):=
\inf_{\phi\in W^{1,p}_{per}(\square,\R^m)}
\int_\square f_{hom}^{[n]}\negthinspace\left(x,y^1, \ldots,y^{n-1}, z+\g \phi(y^{n-1})\right)dy^{n-1},\\ 
&\quad\vdots\\
&f_{hom}\left(x,z\right)= f_{hom}^{[1]}\left(x,z\right):=
\inf_{\phi\in W^{1,p}_{per}(\square,\R^m)}
\int_\square f_{hom}^{[2]}\negthinspace\left(x,y^1, z+\g \phi(y^1)\right)dy^1. 
\end{split}\end{equation*}

\vspace{6pt}
\noindent
\begin{rem}
(i)The convexity, the $p$-coerciveness and the $p$-growth condition on $f$ give the 
corresponding properties for the function $f_{hom}^{[n]}$. 
Moreover, $f_{hom}^{[n]}$ is still an admissible integrand. In fact,
for every $\delta>0$ there exist a compact set $X\subseteq\Om$ with $\snorm{\Om\backslash X}\leq\delta$ 
and a compact set $Y\subseteq\square$ with $\snorm{\square\backslash Y}\leq\delta$,
such that the restriction of $f$ to $X\times Y^{n-1}\times\square\times\M^{m\times d}$ is 
continuous in $(x,y^1,\ldots,y^{n-1},z)$ for \textit{a.e.} $y^n\in\square$.
Consequently, following closely \cite[Lemma 4.1]{Fons03}, it can be proved that
$f_{hom}^{[n]}$ is continuous on  $X\times Y^{n-1}\times\M^{m\times d}$.\\
(ii)Clearly, the properties of $f_{hom}^{[n]}$ give the corresponding ones for $f_{hom}^{[n-1]}$ and so on.
\end{rem}

\vspace{6pt}
We prove only the inequality $f_{hom}\leq f_{hom}^{[1]}$,
since the opposite inequality comes directly.
Fixed $(x,z)\in\Om\times\M^{m\times d}$ and 
$\phi_k\in\Phi_{k,reg}$ for $k=1,\ldots,n-1$,
by using a commutation argument as Lemma \ref{lemma commutativita}, we get
\begin{flalign*}
\inf_{\phi_n\in\Phi_{n,reg}}&\int_{\square^n}
f\Bigl(x, y, z+\sum_{k=1}^n\g_{y^k}\phi_k(\yk)\Bigr)\,dy\\
=&\int_{\square^{n-1}}\inf_{\phi_n\in\Phi_{n,reg}}
\left(\int_\square f\Bigl(x, y, z+\sum_{k=1}^n\g_{y^k}\phi_k(\yk)\Bigr)\,dy^n\right)dy^1\ldots dy^{n-1}\\
\leq&\int_{\square^{n-1}}
f_{hom}^{[n]}\Bigl(x,y^1, \ldots,y^{n-1}, z+\sum_{k=1}^{n-1}\g_{y^k}\phi_k(\yk)\Bigr)\,dy^1\ldots dy^{n-1}\,.
\end{flalign*}

\vspace{6pt}
\noindent
By repeating the commutation procedure, we obtain
\begin{flalign*}
\inf_{\phi_{n-1}\in\Phi_{n-1,reg}}&\int_{\square^{n-1}}
f_{hom}^{[n]}\Bigl(x,y^1, \ldots,y^{n-1}, z+\sum_{k=1}^{n-1}\g_{y^k}\phi_k(\yk)\Bigr)\,dy^1\ldots dy^{n-1}\\
\leq&\int_{\square^{n-2}}
f_{hom}^{[n-1]}\Bigl(x,y^1, \ldots,y^{n-2}, z+\sum_{k=1}^{n-2}\g_{y^k}\phi_k(\yk)\Bigr)\,dy^1\ldots dy^{n-2}
\end{flalign*}
and so on. Then 
\begin{equation*}
f_{hom}(x,z)\leq\inf_{\phi_1\in\Phi_{1,reg}}\ldots\inf_{\phi_n\in\Phi_{n,reg}}
\int_{\square^n}f\Bigl(x, y, z+\sum_{k=1}^n\g_{y^k}\phi_k(\yk)\Bigr)\,dy\leq f_{hom}^{[1]}(x,z)\,.
\end{equation*}

\bigskip
\bigskip
\centerline{\textsc{Acknowledgments}}
\bigskip
\noindent
I wish to thank Gianni Dal Maso for many helpful and interesting discussions.
\bigskip


\end{document}